\documentclass[a4paper, 11pt]{article}

\usepackage{amsmath,amssymb,amsthm}
\usepackage[latin1]{inputenc}
\usepackage{version,tabularx,multicol}
\usepackage{graphicx,float,psfrag}
\usepackage{stmaryrd}
 \usepackage{enumerate}
 \usepackage{color}
\usepackage[pdftex]{hyperref}
\usepackage{authblk}

\usepackage[numbers,sort&compress]{natbib}

\theoremstyle{plain}
\newtheorem{theorem}{Theorem}[section]

 \newtheorem{lemma}[theorem]{Lemma}

\newcommand{\bn}[1]{\textbf{#1}}
\newcommand{\bm}[1]{\emph{\textbf{#1}}}

  \makeatletter
  \@addtoreset{equation}{section}
  \makeatother

 \def\beqlb{\begin{eqnarray}}\def\eeqlb{\end{eqnarray}}
 \def\beqnn{\begin{eqnarray*}}\def\eeqnn{\end{eqnarray*}}

 \def\qed{\hfill$\Box$\medskip}

\newcommand{\bcen}{\begin{center}}
\newcommand{\ecen}{\end{center}}
\newcommand{\bgeqn}{\begin{equation}}
\newcommand{\edeqn}{\end{equation}}

\begin{document}

\title{Characterization of fixed points of infinite-dimensional generating functions}

\author{ Jiangrui Tan\thanks{%
School of Mathematical Sciences \& Laboratory of Mathematics and Complex
Systems, Beijing Normal University, Beijing 100875, P.R. China. Email:
jr\_tan@mail.bnu.edu.cn} \ and Mei Zhang\thanks{%
Corresponding
author. School of Mathematical Sciences \& Laboratory of Mathematics and Complex
Systems, Beijing Normal University, Beijing 100875, P.R. China. Email:
meizhang@bnu.edu.cn} }
\maketitle

\noindent{\bf Abstract}\quad
 This paper is concerned with the characterizations of fixed points of the generating function of branching processes with countably infinitely many types. We assume each particle of type $i$ can only give offspring of type $j\geq i$,  whose number only depends on $j-i$. We prove that, for these processes, there are at least countably infinitely many fixed points of the offspring generating function, while the extinction probability set of the process has only $2$ elements. This phenomenon contrasts sharply with those of finite-type branching processes. Our result takes one step forward on the related conjecture on the fixed points of infinite-dimensional generating functions in literature. In addition, the asymptotic behavior of the components of fixed point is given.
 \vspace{0.3cm}



\bigskip

\section{ Introduction and Preliminaries}
\vspace{3mm}

 \ \ Galton-Watson branching processes (GWBPs) are models describing the evolution of particle systems where independent particles reproduce and die. If the reproduction law varies in some classes of particles, multi-type GWBPs are suitable models (see more details in \cite[Chapter 5]{Athreya1972Branching}). In this paper, we focus on the GWBPs with countably many types, which can naturally be interpreted as branching random walks on an infinite graph where the types of particles correspond to the vertices of graph (see \cite{SZ}). These processes are of many applications, especially used as stochastic models for biological populations (see \cite{KA}).

  We consider a GWBP with countably many types $\{\bm{Z}_n;n\geq0\}$ in which the generating function $\bm{F}(\bm{s})=(F^{(1)}(\bm{s}),F^{(2)}(\bm{s}),\cdots)$ has the form as
\begin{equation}\label{fs}
F^{(i)}(\bm{s})=\sum_{j_1,j_2,\cdots\geq0}P(j_1,j_2,\cdots)\prod_{k=1}^{\infty}s_{i+k-1}^{j_k},
\end{equation}
where $\bm{s}=(s_1,s_2,\cdots)$ and $P(j_1,j_2,\cdots)$ represents the probability of a particle of type $i$ gives $j_k$ offspring of type $i+k-1$ for $k\geq1$ respectively. As we see, the reproduction law only depends on the value of variance of types (or says, the distance between two vertices).

  Let $\bn{1}$ and $\bn{0}$ be the infinite vectors of $1$s and $0$s. Denote the mean matrix of $\{\bm{Z}_n;n\geq0\}$ by $\bm{M}=((m_{ik}))$ where
\[
m_{ik}=\cfrac{\partial\,F^{(i)}}{\partial\,s_k}(\bn{1}).
\]
Clearly $m_{ik}=0$ if $k<i$. For $k\geq1$, define
\begin{equation}\label{mk}
M_k=\frac{\partial F^{(1)}}{\partial s_k}(\bn{1})=m_{1k} ~~\mbox{and}~~ M=\sum_{k\geq1}M_k.
\end{equation}
 To avoid trivialities, we make the following basic assumptions:
 \begin{flalign*}
&\emph{\bf \text{A1:}}~\text{For any~} k\geq i\ge 1, \mbox{there exists a positive integer~} n ~\mbox{such that~} (\bm{M}^n)_{ik}>0.&\\
&\emph{\bf \text{A2:}}~P(\bn{0})>0~\text{and}~\mathbb{P}(|\bm{Z}_1|>1)>0.&\\
&\emph{\bf \text{A3:}}~ M_1<1~\text{and}~M<\infty.&
\end{flalign*}\par
For a general GWBP with countably infinitely many types $\{\bm{X}_n;n\geq0\}$ with offspring generating function $\emph{\textbf{g}}(\emph{\textbf{s}})$, the extinction can be of the whole population--global extinction, in all finite subsets of types--partial extinction, or more generally,
in any fixed subset of types $A$--local extinction in $A$. More precisely, let $\mathcal{T}\subset \mathbb{N}=\{1,2,\cdots\}$. Then the local extinction probability $\bm{q}(\mathcal{T})=\{q^{(i)}(\mathcal{T});i\geq1\}$ in $\mathcal{T}$ is defined as
\begin{eqnarray*}
q^{(i)}(\mathcal{T})=\mathbb{P}\left(\lim_{n\rightarrow\infty}\sum_{l\in\mathcal{T}}X^{(l)}_n=0\bigg|\bm{X}_0=\bm{e}_i\right),
\end{eqnarray*}
where $X_n^{(l)}$ is the $l$-th component of $\bm{X}_n$, $\bm{e}_i$ is the infinite vector where all entries equal to zero except that entry $i$ equals to $1$.

$\bm{q}(\mathcal{T})$ is called global extinction probability when $\mathcal{T}=\mathbb{N}$, and called partial extinction probability when $\mathcal{T}$ is a finite set. In irreducible cases, $\bm{q}(\mathcal{T})$ coincides for any finite subset $\mathcal{T}$ (see more details for $\bm{q}(\mathcal{T})$ in \cite{GF}).
For every $\mathcal{T}\subset\mathbb{N}$, it is easy to know $\bm{q}(\mathcal{T})$ is the solution of $\emph{\textbf{g}}(\emph{\textbf{s}})=\emph{\textbf{s}}$. If we denote the extinction probability set and fixed point set of $\emph{\textbf{g}}(\cdot)$ by
$$
\Theta=\{\bm{q}(\mathcal{T}):\mathcal{T}\subset\mathbb{N}\}~\mbox{and}
~\Lambda=\{\bm{s}\in[0,1]^{\mathbb{N}}:\bm{g}(\bm{s})=\bm{s}\},
$$
respectively. It is clear that $\Theta\subset\Lambda$. The characterizations of $\Theta$ and $\Lambda$ are of independent interest.
\par
It is known that for irreducible finite-type branching processes, $\Theta$ and $\Lambda$ are well established. That is, $\Theta=\Lambda=\{\bm{q},\bn{1}\}$ where $\bm{q}$ is the extinction probability (see \cite{Athreya1972Branching}). When
the set of type is countably infinite, the characterizations of $\Theta$ and $\Lambda$ become more complicated. Moyal~\cite{J1962Multiplicative} shows that the global extinction probability $\bm{q}(\mathbb{N})$ is the minimal element in $\Lambda$. Bertacchi et al~\cite{GF} show that in irreducible cases, the partial extinction probability $\tilde{\bm{q}}$ is either the maximal element of $\Lambda$ (equals to $\bn{1}$) or the second large element of $\Lambda$ (see \cite[Theorem 3.1]{GF}). Braunsteins and Hautphenne~\cite{2017arXiv170602919B} prove that, for a class
of branching processes with countably many types called lower Hessenberg branching
processes (i.e., a type $i$ particle can only give  offspring of type $j\leq i+1$), there exists a continuum of fixed points between $\bm{q}(\mathbb{N})$ and $\tilde{\bm{q}}$ if $\bm{q}(\mathbb{N})<\tilde{\bm{q}}\leq\bn{1}$.

In section 4 of \cite{BZ}, Bertacchi and Zucca make a detailed
summary for the known  results related to $\tilde{\bm{q}}$, $\Lambda$ and $\Theta$, and list some open questions on the characterizations of $\Lambda$ and $\Theta$. In particular, they make a conjecture that $\Lambda$ (or $\Theta$) is either finite or uncountable, which is also raised  similarly by \cite[Conjecture 5.1]{BH}. In this paper, we make a step further on the problem by proving that, for an infinite-type GWBP with the offspring generating function of the form (\ref{fs}), if $M\leq1$, then $\Theta=\Lambda=\{\bn{1}\}$. If $M>1$,
then $\Lambda$ has at least countably many fixed points while $\Theta=\{q\bn{1},\bn{1}\}$.

Before stating our main result, we make a brief discussion for the extinction of $\{\bm{Z}_n;n\geq0\}$ with generating function (\ref{fs}). Due to $M_1<1$, the type of descendants of any type $i~(i\ge 1)$ particle will exceed $i$ in finite time which implies the extinction in any finite typeset $\mathcal{T}$. Hence the partial extinction probability $\tilde{\bm{q}}=\bn{1}$. On the other hand, if we ignore the type of each particle, then $\{\bm{Z}_n;n\geq0\}$ degenerates to a classical GWBP with offspring p.g.f.
\[
F_0(s)=\sum_{k=0}^{\infty}\sum_{|\bm{j}|=k}P(\bm{j})s^{k}.
\]
Clearly if and only if $M=F_0'(1)>1$, there exists a unique solution in $(0,1)$ to the equation $F_0(s)=s$ which we denote by $q$, where $q>0$ follows by $P(\bn{0})>0$. Since the offspring distribution only relies on the variance of types, the global extinction probability $\bm{q}(\mathbb{N})=q\bn{1}$. Next, let $\mathcal{T}^\ast$ be an arbitrary infinite typeset.  By assumption $\bn{A1}$, it is not difficult to see $\bm{q}(\mathcal{T}^\ast)= q\bn{1}$. Hence $\Theta=\{\bn{1}\}$ if $M\leq1$, and $\Theta=\{q\bn{1},\bn{1}\}$ if $M>1$.

Now we state the main result of this paper.

\begin{theorem}\label{th1}
If $M\leq1$, then $\Theta=\Lambda=\{\bm{1}\}$. If $M>1$,
then $\Lambda$ has at least countably infinitely many fixed points while $\Theta=\{q\bm{1},\bm{1}\}$. Moreover, for any  $\bn{r}=(r^{(i)})_{i\ge 1}\in\Lambda\setminus\Theta$, it holds that
\[
\lim_{i\to\infty}\frac{1-r^{(i+1)}}{1-r^{(i)}}=\gamma,
\]
where $\gamma$ is the unique solution in $(0,1)$ to the equation $\sum_{i=1}^{\infty}M_is^{i-1}=1$.
\end{theorem}

\section{Proofs}

At first, we make a paraphrasing for $F^{(1)}(\bm{s})$ which is of many uses in this paper.\par
 Let $h_0=P(\bn{0})$ and $h_i~(i\geq1)$ be the probability that the maximal index of the offspring type  of a type $1$ particle is $i$. That is,
\[
h_i=\sum_{\begin{subarray}{c}j_1,\cdots,j_{i-1}\geq0\\j_{i}>0\end{subarray}}P(j_{1},\cdots,j_{i},0,0,\cdots).
\]
Define~the $k$-dimensional probability generating function
\[
f_k(s_1,\cdots,s_k)=\sum_{\begin{subarray}{c}j_1,\cdots,j_{k-1}\geq0\\j_k>0\end{subarray}}\frac{P(j_1,\cdots,j_k,0,\cdots)}
{h_k}\prod_{i=1}^{k}s_i^{j_i}.
\]
Then
\begin{equation}\label{eq}
F^{(1)}(\bm{s})=h_0+\sum_{k=1}^{\infty}h_kf_k(s_1,\cdots,s_k).
\end{equation}
For $k\geq j\geq 1$, define
\[
a_{k,j}=\frac{\partial f_k}{\partial s_j}(\bn{1}).
\]
 Then calculation yields
\[
M_i=\sum_{k=i}^{\infty}h_ka_{k,i}.
\]

Here are some notations for simplicity. Throughout this paper, we use bold characters to represent vectors or matrices and $x^{(i)}$ to denote the $i$-th component of any vector $\bm{x}\in[0,1]^{\mathbb{N}}$. For $j\geq i$, define
\[
\bm{x}_{i\rightarrow j}=(x^{(i)},x^{(i+1)},\cdots,x^{(j)}).
\]
For any $s_1,s_2,\cdots,s_k~(k>0)$, define $(s_1,s_2,\cdots,s_k,\bm{x})=(s_1,s_2,\cdots,s_k,x^{(1)},x^{(2)},\cdots)$ and write $F^{(1)}(s_1,s_2,\cdots,s_k,\bm{x})=F^{(1)}\big((s_1,s_2,\cdots,s_k,\bm{x})\big)$ (Similarly for other  vector functions). Write $\bm{x}\leq\bm{y}$ ($\bm{x}<\bm{y}$) if $x^{(i)}\leq y^{(i)}$ ($x^{(i)}<y^{(i)}$) for all $i\ge 1$. In addition, we write $\bn{1}-\bm{x}_{i\rightarrow j}=(1-x^{(i)},\cdots,1-x^{(j)})$ for any $j\geq i\ge 1$ and $\bm{x}^T$ as the transpose of $\bm{x}$.
\par
It is known that for any multi-type GWBP with generating function $\bm{G}(\cdot)$ and mean matrix $\bm{M}_0$, it holds that
\begin{equation}\label{ex}
(\bn{1}-\bm{G}(\bm{s}))^T=(\bm{M}_0-\bm{E}(\bm{s}))(\bn{1}-\bm{s})^T,
\end{equation}
where $\bn{0}\leq\bm{E}(\bm{s})\leq\bm{M}_0$ elementwise, $\bm{E}(\bm{s})$ is non-increasing in $\bm{s}$ (with respect to the partial order induced by ``$\leq$") and tends to $\bn{0}$ as $\bm{s}\rightarrow\bn{1}$ (see the proof of Theorem 1 on \cite[Page 414]{Joffe}).
By (\ref{ex}), for any $k>0$, we have
\[
1-f_k(\bn{1}-\bm{s}_{1\rightarrow k})=\sum_{i=1}^{k}\left(\frac{\partial f_k}{\partial s_i}(\bn{1})-E_{k,i}(\bn{1}-\bm{s}_{1\rightarrow k})\right)s_i=\sum_{i=1}^{k}(a_{k,i}-E_{k,i}(\bn{1}-\bm{s}_{1\rightarrow k}))s_i
\]
holds for some $\{E_{k,i}(\cdot);k\geq i\geq1\}$, where
\begin{eqnarray*}
E_{k,i}(\bn{1}-\bm{s}_{1\rightarrow k})&=&\sum_{\begin{subarray}{c}j_1,\cdots,j_{k-1}\geq0\\j_k>0\end{subarray}}\frac{P(j_1,\cdots,j_k,0,\cdots)}
{h_k}j_i\left[1-\int_{0}^1\frac{\prod_{l=1}^{k}(1-s_{l}x)^{j_l}}{1-s_{i}x}dx\right]
\end{eqnarray*}
by applying equation (4.3) of \cite[Page 414]{Joffe}.

The following lemma shows that for any given $\bm{y}\in(0,1)^{\mathbb{N}}$, we can find $x\in(0,1)$ such that $F^{(1)}(x,\bm{y})=x$.
\begin{lemma}\label{lem1}
Given $\textbf{y}\in(0,1)^{\mathbb{N}}$, there is a unique solution in $(0,1)$ to the equation $F^{(1)}(x,\textbf{y})=x$. Denote the solution by $F_{-1}^{(1)}(\textbf{y})$. Moreover, if $\textbf{y}_1\geq\textbf{y}_{2}$, then $F_{-1}^{(1)}(\textbf{y}_1)\geq F_{-1}^{(1)}(\textbf{y}_2)$ and further if $y_1^{(1)}>y_2^{(1)}$, then $F_{-1}^{(1)}(\textbf{y}_1)> F_{-1}^{(1)}(\textbf{y}_2)$.
\end{lemma}

\proof
Given $\bm{y}\in(0,1)^{\mathbb{N}}$, by assumption $\bn{A1}$, we have $F^{(1)}(0,\bm{y})>0$ and $F^{(1)}(1,\bm{y})<1$, and then the first part of the lemma follows.\par

Noticing that if $\bm{y}_1\geq\bm{y}_{2}$, then $F^{(1)}(x,\bm{y}_1)\geq F^{(1)}(x,\bm{y}_2)$ for any $x\in(0,1)$. Since $F^{(1)}(x,\bm{y})$ is continuous and increasing with respect to $x$ for any given $\bm{y}\in(0,1)^{\mathbb{N}}$, we have  $F_{-1}^{(1)}(\bm{y}_1)\geq F_{-1}^{(1)}(\bm{y}_2)$. From assumption $\bn{A1}$, a type $i$ particle has a positive probability to give type $i+1$ particles in one generation. Hence if $y_1^{(1)}>y_2^{(1)}$, $F^{(1)}(x,\bm{y}_1)> F^{(1)}(x,\bm{y}_2)$ for any $x\in(0,1)$ which implies $F_{-1}^{(1)}(\bm{y}_1)> F_{-1}^{(1)}(\bm{y}_2)$ and the lemma follows.\qed

\begin{lemma}\label{lem2}
If $M>1$, there is a unique solution in $(0,1)$ to the equation $G(s):=\sum_{i=1}^{\infty}M_is^{i-1}=1$.
\end{lemma}
\proof
The lemma follows from $1<G(1)<\infty$, $G(0)<1$ and $G(s)$ is continuous in $(0,1)$.\qed

Denote the unique solution of $G(s)=1$ in $(0,1)$ by $\gamma$. Define the vector set
\[
\mathcal{H}(\gamma)=\left\{\bm{x}:\forall~i\ge 1,x^{(i)}\in (0,1), x^{(i)}>x^{(i+1)}~\mbox{and}~\lim_{i\rightarrow\infty}\frac{x^{(i+1)}}{x^{(i)}}=\gamma\right\}.
\]
Obviously $\mathcal{H}(\gamma)\subset l_2$, where~$l_2$ is the normalized sequence space $\{\bm{x}:\sum_{i\ge 1}|x^{(i)}|^2<\infty\}$. For any $\bm{x}\in l_2$, define
\[
\bm{T}(\bm{x})=\bn{1}-\bm{F}(\bn{1}-\bm{x}).
\]

\begin{lemma}\label{lem3}
For any $\bn{x}\in\mathcal{H}(\gamma)$, $\bn{T}(\bn{x})\in\mathcal{H}(\gamma)$.
\end{lemma}
\proof
From (\ref{eq}) and (\ref{ex}),
\begin{eqnarray}
\frac{T^{(i)}(\bm{x})}{x^{(i)}}&=&\frac{1-F^{(1)}(\bn{1}-\bm{x}_{i\rightarrow\infty})}{x^{(i)}}\nonumber\\
&=&\frac{\sum_{k=1}^{\infty}h_k\big(1-f_k(\bn{1}-\bm{x}_{i\rightarrow i+k-1})\big)}{x^{(i)}}\nonumber\\
&=&\sum_{k=1}^{\infty}h_k\sum_{j=1}^{k}\big(a_{k,j}-E_{k,j}(\bn{1}-\bm{x}_{i\rightarrow i+k-1})\big)\frac{x^{(i+j-1)}}{x^{(i)}}\nonumber\\
&\leq&\sum_{k=1}^{\infty}h_k\sum_{j=1}^{k}a_{k,j}=M<\infty.\label{eq2}
\end{eqnarray}
 Then from the dominated convergence theorem,
\begin{eqnarray*}
\lim_{i\rightarrow\infty}\frac{T^{(i)}(\bm{x})}{x^{(i)}}
&=&\lim_{i\rightarrow\infty}\left(h_1a_{1,1}+h_2\left(a_{2,1}+a_{2,2}\gamma \right)+ h_3\left(a_{3,1}+a_{3,2}\gamma+a_{3,3}\gamma^2\right)+\cdots\right)\\
&=& M_1+M_2\gamma+M_3\gamma^2+\cdots\\
&=&1.
\end{eqnarray*}
Hence
\[
\lim_{i\rightarrow\infty}\frac{T^{(i+1)}(\bm{x})}{T^{(i)}(\bm{x})}=
\lim_{i\rightarrow\infty}\frac{T^{(i+1)}(\bm{x})}{x^{(i+1)}}\cdot\frac{x^{(i)}}{T^{(i)}(\bm{x})}\cdot\frac{x^{(i+1)}}{x^{(i)}}=\gamma.
\]
Noting that $T^{(i+1)}(\bm{x})=T^{(i)}(\bm{x}_{2\rightarrow\infty})$, then $T^{(i)}(\bm{x})>T^{(i+1)}(\bm{x})$ follows by $\bm{x}>\bm{x}_{2\rightarrow\infty}$ and the lemma follows.\qed\par
The following lemma shows $T^{(i)}(\cdot)$ is pointwisely continuous for any $i\ge 1$.

\begin{lemma}\label{lem4}
For any $\textbf{x}_1,\textbf{x}_2\in l_2$, it holds that for any $i\ge 1$,
\[
|T^{(i)}(\textbf{x}_1)-T^{(i)}(\textbf{x}_2)|\leq M\cdot|x_1^{(i)}-x_2^{(i)}|.
\]
\end{lemma}
\proof
  From \cite[Page 57, Section 3]{BZC}, $\bm{F}(\cdot)$ is continuous with respect to pointwise convergence topology. Hence $\bm{T}(\cdot)$ is pointwisely continuous and the lemma follows. \qed

From now on, we fix $\bm{x}$ for some $\bm{x}\in\mathcal{H}(\gamma)$. Define
\begin{eqnarray}
& &\eta_n^{[1]}=\eta^{[1]}_{n}(\bm{x})=
T_{-1}^{(1)}(\bm{x}_{n\rightarrow\infty}):=1-F^{(1)}_{-1}(\bn{1}-\bm{x}_{n\rightarrow\infty})~~\mbox{and}\nonumber\\
& &
\eta_n^{[i]}=T_{-1}^{(1)}(\eta_n^{[i-1]},\eta_n^{[i-2]},\cdots,\eta_n^{[1]},\bm{x}_{n\rightarrow\infty}) ~~\mbox{for}~i\geq2.\label{eta1}
\end{eqnarray}
 Observing that by Lemma~\ref{lem1}, we construct the sequence  $\{\eta_n^{[k]};1\le k \le i\}$  such that
\begin{equation}\label{k11}
(\eta_n^{[i]},\eta_n^{[i-1]},\cdots,\eta_n^{[1]},\bm{x}_{n\rightarrow\infty})=T^{(k)}(\eta_n^{[i]},\eta_n^{[i-1]},\cdots,\eta_n^{[1]},\bm{x}_{n\rightarrow\infty})
\end{equation}
hold for all $1\le k \le i$.

 The following lemma shows that for $n$ large enough, $\eta_n^{[i]}$ is strictly increasing with respect to $i$.
\begin{lemma}\label{lem6}
There exists positive integer $N_0$, such that for $n\geq N_0$, we have $\eta_n^{[1]}>x^{(n)}$ and $\eta_n^{[i]}>\eta_{n}^{[i-1]}$ for all $i\geq2$.
\end{lemma}
\proof
Note that  by (\ref{k11}),
\begin{eqnarray*}
\frac{\eta_n^{[1]}}{x^{(n)}}&=&\frac{1-F^{(1)}(1-\eta_n^{[1]},\bn{1}-\bm{x}_{n\rightarrow\infty})}{x^{(n)}}\\
&=&\frac{h_1\big(1-f_1(1-\eta_n^{[1]})\big)+\sum_{k=2}^{\infty}h_k\big(1-f_k(1-\eta_n^{[1]},\bn{1}-\bm{x}_{n\rightarrow n+k-2})\big)}{x^{(n)}}\\
&=&\sum_{k=1}^{\infty}h_k\bigg[\big(a_{k,1}-E_{k,1}(1-\eta_n^{[1]},\bn{1}-\bm{x}_{n\rightarrow n+k-2})\big)\frac{\eta_n^{[1]}}{x^{(n)}}+\\
& &\sum_{j=2}^{k}\big(a_{k,j}-E_{k,j}(1-\eta_n^{[1]},\bn{1}-
\bm{x}_{n\rightarrow n+k-2})\big)\frac{x^{(n+j-2)}}{x^{(n)}}\bigg],
\end{eqnarray*}
provided $E_{1,1}(1-\eta_n^{[1]},\bn{1}-\bm{x}_{n\rightarrow n-1})=E_{1,1}(1-\eta_n^{[1]})$ and $\sum_{j=2}^{1}=0$. Then
\begin{eqnarray}
\frac{\eta_n^{[1]}}{x^{(n)}}&=&\left(\sum_{k=2}^{\infty}h_k\sum_{j=2}^{k}\big(a_{k,j}-E_{k,j}(1-\eta_n^{[1]},\bn{1}-\bm{x}_{n\rightarrow n+k-2})\big)\frac{x^{(n+j-2)}}{x^{(n)}}\right)\nonumber\\
& &\boldsymbol{\cdot}\left(1-\sum_{k=1}^{\infty}h_k\big(a_{k,1}-E_{k,1}(1-\eta_n^{[1]},\bn{1}-\bm{x}_{n\rightarrow n+k-2})\big)\right)^{-1}.\label{eq1}
\end{eqnarray}
Noting that $x^{(i)}>x^{(i+1)}$ for all $i\ge 1$, we have
 \begin{eqnarray*}
& & 0<\sum_{k=2}^{\infty}h_k\sum_{j=2}^{k}[a_{k,j}-E_{k,j}(1-\eta_n^{[1]},\bn{1}-\bm{x}_{n\rightarrow n+k-2})]\frac{x^{(n+j-2)}}{x^{(n)}}\le \sum_{k=2}^{\infty}h_k\sum_{j=2}^{k}a_{k,j}<\infty,\\
& & 0<\sum_{k=1}^{\infty}h_k[a_{k,1}-E_{k,1}(1-\eta_n^{[1]},\bn{1}-\bm{x}_{n\rightarrow n+k-2})] \le \sum_{k=1}^{\infty}h_ka_{k,1}<1.
\end{eqnarray*}
 Observe that $\eta_n^{[1]}\to 0$, and for each $k$, $\bm{x}_{n\rightarrow n+k-2}\to \bn{0}$ as $n\to \infty$. Then by the dominated convergence theorem,
\begin{eqnarray*}
\lim_{n\rightarrow\infty}\frac{\eta_n^{[1]}}{x^{(n)}}=\left(\sum_{k=2}^{\infty}h_k\sum_{j=2}^{k}a_{k,j}\gamma^{j-2}\right)\left(1-\sum_{k=1}^{\infty}h_ka_{k,1}\right)^{-1}
=\gamma^{-1},
\end{eqnarray*}
where the last equality follows by $$G(\gamma)=\sum_{k=1}^{\infty}h_k\sum_{j=1}^{k}a_{k,j}\gamma^{j-1}=1.$$
 Since $\gamma<1$, there exists $N_0$ such that $\eta_n^{[1]}>x^{(n)}$ for $n\geq N_0$. Next, from Lemma \ref{lem1} and (\ref{eta1}), $\eta_n^{[2]}>\eta_n^{[1]}$ and $\eta_n^{[i]}>\eta_n^{[i-1]}$ follows consequently for $i\geq2$. The proof is completed.\qed

Define the set of mapping on positive integers:
\[
\mathcal{A}=\left\{I:\mathbb{N}\mapsto\mathbb{N}, I(n+1)>I(n)~\mbox{for}~ n\geq1\right\}.
\]
For any $I, J\in\mathcal{A}$, clearly $\lim_{n\rightarrow\infty}I(n)=\lim_{n\rightarrow\infty}J(n)=+\infty$. We can define a sequence of vectors $\{\bm{y}_n(I,J);n\geq1\}$, where
\begin{equation}\label{ynij}
\bm{y}_n(I,J)=(\eta_{J(n)}^{[I(n)]},\eta_{J(n)}^{[I(n)-1]},\cdots,\eta_{J(n)}^{[1]},x^{(J(n))},x^{(J(n)+1)},\cdots).
\end{equation}

 From (\ref{k11}), $y_n^{(k)}(I,J)=T^{(k)}(\bm{y}_n(I,J))$ for all $1\le k\le I(n)$. Hence, from now on, we will prove that there exist $I, J\in\mathcal{A}$ such that $\bm{y}_n(I,J)$ converges to some fix point of~$\bm{T}(\cdot)$ pointwisely. At first, we have the following lemma:
\begin{lemma}\label{lem7}
For any~$I,~J\in\mathcal{A}$,  ~$\lim_{n\rightarrow\infty}\parallel\bn{y}_n(I,J)-\bn{T}\big(\bn{y}_n(I,J)\big)\parallel_{l_2}=0$.
\end{lemma}
\proof
 From the  definition of $\eta_n^{[i]}$ and (\ref{ynij}),  for any~$1\leq i\leq I(n)$,
\begin{eqnarray*}
y_n^{(i)}(I,J)&=&\eta_{J(n)}^{[I(n)-i+1]}=T^{(1)}(\eta_{J(n)}^{[I(n)-i+1]},\eta_{J(n)}^{[I(n)-i]},\cdots,\eta_{J(n)}^{[1]},x^{(J(n))},\cdots)\\
&=&T^{(i)}\big(\bm{y}_n(I,J)\big).
\end{eqnarray*}
 For $i\ge I(n)+1$, also by (\ref{ynij}), we have $$y_n^{(i)}(I,J)=x^{(i-I(n)-1+J(n))}, \ \ T^{(i)}\big(\bm{y}_n(I,J)\big)=T^{(1)}(\bm{x}_{(i-I(n)-1+J(n))\rightarrow\infty}). $$ Therefore,
\begin{eqnarray}\label{ynij1}
\parallel\bm{y}_n(I,J)-\bm{T}\big(\bm{y}_n(I,J)\big)\parallel_{l_2}^2&=& \sum_{k=J(n)}^\infty\big[ x^{(k)}-T^{(1)}(\bm{x}_{k\rightarrow\infty})\big]^2.
\end{eqnarray}
From Lemma \ref{lem3},
\[
\sum_{k=0}^{\infty}\big(x^{(k)}-T^{(1)}(\bm{x}_{k\rightarrow\infty})\big)^2<\infty.
\]
Letting $J(n)\to \infty$ in (\ref{ynij1}),
we complete the proof.\qed

  In the following, Lemmas \ref{lem8} and \ref{lem10} show that there exist  $\bm{y}\in (0,1]^\mathbb{N}$, $I_1, J_1\in\mathcal{A}$, such that $\bm{y}\neq(1-q)\bn{1}$, and $\bm{y}_n(I_1,J_1)$ converges to $\bm{y}$ pointwisely.
\begin{lemma}\label{lem8}
There exist $I_0, J_0\in \mathcal{A}$, such that $\lim_{n\rightarrow\infty}\parallel\bn{y}_n(I_0,J_0)\parallel_{l_2}$ exists and $\in(0,\infty)$.
\end{lemma}
\proof
For~$j\geq N_0$, define
\begin{equation}\label{yij}
\bm{y}[i,j]=(\eta_j^{[i]},\eta^{[i-1]}_{j},\cdots,\eta_j^{[1]},x^{(j)},x^{(j+1)},\cdots).
\end{equation}
From Lemma \ref{lem1} and the definition of $\eta_n^{[i]}$, we have $\eta_{j+1}^{[1]}<\eta_{j}^{[1]}$ and hence $\eta_{j+1}^{[i]}<\eta_j^{[i]}$ for any $i>1$. Then $\parallel\bm{y}[i,j]\parallel_{l_2}$ is strictly decreasing to $0$ with respect to~$j$ (for fixed $i$). On the other hand, it is obvious that $\parallel\bm{y}[i,j]\parallel_{l_2}$ is strictly increasing to $+\infty$ with respect to $i$ (for fixed $j$).\par
Let $y_0:=\parallel\bm{y}[0,N_0]\parallel_{l_2}=\parallel\bm{x}_{N_0\rightarrow\infty}\parallel_{l_2}$. Then there exists~$m_1>N_0$, such that $\parallel\bm{y}[1,m_1]\parallel_{l_2}<y_0$. Next, by Lemma~\ref{lem6}, we have $\parallel\bm{y}[i,m_1]\parallel_{l_2}$ is increasing to $+\infty$ with respect to $i$. Then there exists $k_1> 1$ such that
\[
\parallel\bm{y}[k_1,m_1]\parallel_{l_2}>y_0~\mbox{while}~
\parallel\bm{y}[k_1-1,m_1]\parallel_{l_2}\leq y_0.
\]
 Noticing that, for the fixed $k_1$, $\parallel\bm{y}[k_1,n]\parallel_{l_2}\to 0$ as $n\to \infty$, so we can choose $m_2>m_1$ such that $\parallel\bm{y}[k_1,m_2]\parallel_{l_2}<y_0$.
Similarly, there exists integer~$k_2> k_1$, such that
\[
\parallel\bm{y}[k_2,m_2]\parallel_{l_2}>y_0~\mbox{while}~
\parallel\bm{y}[k_2-1,m_2]\parallel_{l_2}\leq y_0.
\]
Therefore, there exist two  sequences of integers $\{k_n\}$ and $\{m_n\}$ which satisfy $k_n>k_{n-1}$ and~$m_n>m_{n-1}$, such that
\[
\parallel\bm{y}[k_n,m_n]\parallel_{l_2}>y_0~\mbox{while}~
\parallel\bm{y}[k_n-1,m_n]\parallel_{l_2}\leq y_0.
\]
 Then by (\ref{yij}) we have that
\begin{eqnarray*}
\parallel\bm{y}[k_n,m_n]\parallel_{l_2}^2&=&\parallel\bm{y}[k_n-1,m_n]\parallel_{l_2}^2
+\left(\eta_{m_n}^{[k_n]}\right)^2\\
&\le& \parallel\bm{y}[k_n-1,m_n]\parallel_{l_2}^2+1\\
&\le& y_0^2+1\\
&<&(y_0+1)^2.
\end{eqnarray*}
 Taking~$I(n)=k_n$, $J(n)=m_n$, then $I,J\in\mathcal{A}$ and $\parallel\bm{y}_n(I,J)\parallel_{l_2}\in(y_0, y_0+1)$. Hence, there exists a subsequence~$\{r_n\}$ such that $\lim_{n\rightarrow\infty}\parallel\bm{y}_{r_n}(I,J)\parallel_{l_2}$ exists and $\in[y_0, y_0+1]$. Choose~$I_0,~J_0$ satisfy~$I_0(n)=I(r_n)$, $J_0(n)=J(r_n)$, clearly $I_0,~J_0\in\mathcal{A}$ and ~$\lim_{n\rightarrow\infty}\parallel\bm{y}_{n}(I_0,J_0)\parallel_{l_2}$ exists and $\in[y_0, y_0+1]$. The proof is completed.\qed

\begin{lemma}\label{lem9}(\cite[Theorem 3.18]{BH1})
A Banach space X is reflexive if and only if every bounded
sequence has a weakly convergent subsequence.
\end{lemma}

The following lemma is a technical lemma for Lemma \ref{lem10}.

\begin{lemma}\label{lem9.5}
Let $\{\alpha_{n,i};i\geq1,n\geq1\}$ be a sequence satisfying $0<\alpha_{\text{inf}}<\alpha_{n,i}\le1$ for some constant $\alpha_{\text{inf}}<\gamma$  for all $n$ and $i$. Providing $\prod_{1}^0=1$, define $$U_{n,i}=\sum_{k=1}^{\infty}M_k\prod_{l=1}^{k-1}\alpha_{n,i+l-1}.$$
If $U_{n,i}$ converges to $1$ as $n\to\infty$ uniformly for $i$, then $\alpha_{n,i}$ converges to $\gamma$ as $n\to \infty $ uniformly for $i$, where $M_k$ and $\gamma$ are defined in (\ref{mk}) and Lemma \ref{lem2}, respectively.
\end{lemma}

\proof
We prove it by contradiction. If the statement {\em ``$\alpha_{n,i}$ converges to $\gamma$ as $n\to \infty $ uniformly for $i$"} is false, then $\limsup_{n\to\infty}\sup_i\alpha_{n,i}>\gamma$ or $\liminf_{n\rightarrow\infty}\inf_i\alpha_{n,i}<\gamma$.

If $\hat{\alpha}:=\limsup_{n\to\infty}\sup_i\alpha_{n,i}>\gamma$ . Choose integer $k_0$ large enough and $\epsilon$ small enough (depends on $k_0$) satisfying
\begin{itemize}
\item[\bf C1:] $\hat{\alpha}-\gamma>\epsilon\left(2+\cfrac{2+3M}{\inf_{k\leq k_0}M_{k}\alpha^{k-2}_{\text{inf}}}\right)$;\\
\item[\bf C2:]$(\hat{\alpha}-\gamma-2\epsilon)(\sum_{j=2}^\infty M_j\alpha_{\text{inf}}^{j-2})>(\gamma-\alpha_{\text{inf}})\frac{M\gamma^{k_0}}{1-\gamma}+\epsilon$,
\end{itemize}
where $M=\sum_{j\geq1}M_j$. We will make use of these conditions in the later and we mention here that such $k_0$ and $\epsilon$ exist and these conditions also hold for any $\epsilon_0<\epsilon$.\par
Since $U_{n,i}$ converges to $1$ as $n\to\infty$ uniformly for $i$, then there exists $\hat{N}_0$ such that for all $n>\hat{N}_0$ and all $i\ge 1$, $U_{n,i}\in(1-\epsilon,1+\epsilon)$. Hence we have $(U_{n,i}-U_{n,i+1})\in(-2\epsilon,2\epsilon)$ for any $i\ge 1$ and $n>\hat{N}_0$, that is
\begin{equation}\label{9.1}
\sum_{k=2}^\infty M_k(\alpha_{n,i}-\alpha_{n,i+k-1})\prod_{l=1}^{k-2}\alpha_{n,i+l}\in(-2\epsilon,2\epsilon).
\end{equation}
Also, by $\sum_{j\geq1}M_j\gamma^{j-1}=1$, we get
\begin{equation}\label{9.2}
U_{n,i}-1=\sum_{k=2}^\infty M_k\left(\prod_{l=1}^{k-1}\alpha_{n,i+l-1}-\gamma^{k-1}\right)\in(-\epsilon,\epsilon)
\end{equation}
holds for all $n>\hat{N}_0$ and $i\ge 1$.\par
On the other hand, there exists $n_0>\hat{N}_0$ such that $\sup_{i}\alpha_{n_0,i}\in(\hat{\alpha}-\epsilon,\hat{\alpha}+\epsilon)$ and $i_0$ (depends on $n_0$) such that $\alpha_{n_0,i_0}>\sup_{i}\alpha_{n_0,i}-\epsilon>\hat{\alpha}-2\epsilon$. Hence
$\alpha_{n_0,i_0}-\alpha_{n_0,i_0+k}>-3\epsilon$
for all $k>0$.
\par
 By (\ref{9.1}), for any $j>1$, if $\alpha_{n_0,i_0+j-1}<\hat{\alpha}-2\epsilon$, then
\begin{eqnarray}
2\epsilon&>&\sum_{k=2}^\infty M_k(\alpha_{n_0,i_0}-\alpha_{n_0,i_0+k-1})\prod_{l=1}^{k-2}\alpha_{n_0,i_0+l}\nonumber\\
&>&M_{j}(\alpha_{n_0,i_0}-\alpha_{n_0,i_0+j-1})\prod_{l=1}^{j-2}\alpha_{n_0,i_0+l}-3\epsilon\cdot M\nonumber\\
&>&M_{j}(\hat{\alpha}-2\epsilon-\alpha_{n_0,i_0+j-1})\alpha_{\text{inf}}^{j-2}-3\epsilon\cdot M\label{99.4}
\end{eqnarray}
where the second inequality follows by $\alpha_{n_0,i}\le1$ for any $i$ and $M=\sum_{k\geq 1}M_k$, the last inequality follows by $\alpha_{n_0,i_0+j-1}<\hat{\alpha}-2\epsilon$ and $\alpha_{n_0,i}>\alpha_{\text{inf}}$ for all $i$. Then
\begin{equation}\label{9.3}
\alpha_{n_0,i_0+j-1}>\hat{\alpha}-\epsilon(2+\frac{2+3M}{M_j\alpha_{\text{inf}}^{j-2}})
\end{equation}
holds for all $j>1$. If $\alpha_{n_0,i_0+j-1}\ge\hat{\alpha}-2\epsilon$, (\ref{9.3}) is obvious. \par
Next, observe that
\begin{eqnarray}
U_{n_0,i_0}-1&=&\sum_{k=2}^\infty M_k\left(\prod_{l=1}^{k-1}\alpha_{n_0,i_0+l-1}-\gamma^{k-1}\right)\nonumber\\
&=&\sum_{k=2}^\infty M_k\left(\sum_{l=1}^{k-1}(\alpha_{n_0,i_0+l-1}-\gamma)\gamma^{l-1}\prod_{j=l}^{k-2}\alpha_{n_0,i_0+j}\right)
\nonumber\\
&=&\sum_{l=1}^\infty (\alpha_{n_0,i_0+l-1}-\gamma)\gamma^{l-1} \left(\sum_{k=l+1}^{\infty}M_k\prod_{j=l}^{k-2}\alpha_{n_0,i_0+j}\right).\label{99.2}
\end{eqnarray}
By (\ref{9.3}) and condition {\bf C1}, $\alpha_{n_0,i_0+l-1}-\gamma>0$ holds for all $1\le l\leq k_0$. Therefore,
\begin{eqnarray}
&&U_{n_0,i_0}-1\nonumber\\
&>&(\alpha_{n_0,i_0}-\gamma) \left(\sum_{k=2}^{\infty}M_k\prod_{j=1}^{k-2}\alpha_{n_0,i_0+j}\right)+\sum_{l=k_0+1}^\infty (\alpha_{n_0,i_0+l-1}-\gamma)\gamma^{l-1} \left(\sum_{k=l+1}^{\infty}M_k\prod_{j=l}^{k-2}\alpha_{n_0,i_0+j}\right)\nonumber\\
&>&(\hat{\alpha}-2\epsilon-\gamma) \left(\sum_{k=2}^{\infty}M_k\alpha_{\text{inf}}^{k-2}\right)+\sum_{l=k_0+1}^\infty (\alpha_{\text{inf}}-\gamma)\gamma^{l-1}\cdot M\nonumber\\
&>&\epsilon, \label{99.1}
\end{eqnarray}
where the last inequality follows by condition {\bf C2}. Meanwhile, by (\ref{9.2}), $U_{n,i}-1<\epsilon$ holds for all $n>\hat{N}_0$ and $i\ge 1$ which leads to a contradiction.

 If $\overline{\alpha}:=\liminf_{n\to\infty}\inf_i\alpha_{n,i}<\gamma$, the discussion is essentially similar to above. The proof is complete.\qed

\begin{lemma}\label{lem10}
 There exist $\bn{y}\in (0,1]^{\mathbb{N}}$,  $I_1, J_1\in\mathcal{A}$, such that $\bn{y}\neq (1-q)\bm{1}$ and $\bn{y}_n(I_1,J_1)$ converges to~$\bn{y}$ pointwisely.
\end{lemma}

\proof
From Lemma \ref{lem8}, $\{\bm{y}_n(I_0,J_0);n\geq1\}$ is a bounded sequence in $l_2$ space. Hence by Lemma \ref{lem9}, there exist a subsequence $\{k_n\}$ and $\bm{y}\in[0,1]^{\mathbb{N}}$ such that $\bm{y}_{k_n}(I_0,J_0)$ converges to $\bm{y}$ pointwisely. Taking $I_1, J_1$ satisfy $I_1(n)=I_0(k_n)$, $J_1(n)=J_0(k_n)$, then $I_1, J_1\in\mathcal{A}$ and $\bm{y}_n(I_1,J_1)$ converges to~$\bm{y}$ pointwisely. Clearly $\bm{y}\neq (1-q)\bn{1}$ follows from $\bm{y}_n(I_0,J_0)$ is bounded in $l_2$. Hence we only need to prove~$\bm{y}\neq\bn{0}$.

$\bullet$ We first prove that, if $y^{(i)}=0$ for some $i$, then $\bm{y}=\bn{0}$.

On the one hand,  by the definition of $\bm{y}_n(I_1,J_1)$, $\eta_n^{[i]}$ and (\ref{k11}), for $1\le i\leq I_1(n)$, we have
\begin{eqnarray}
y_n^{(i)}(I_1,J_1)&=&\sum_{k=1}^{\infty}h_k\sum_{j=1}^{k}\big(a_{k,j}-E_{k,j}(\bn{1}-\bm{y}_n(I_1,J_1))_{i\rightarrow i+k-1}\big)y_n^{(i+j-1)}(I_1,J_1).\label{yni}
\end{eqnarray}
Hence by similar calculation with (\ref{eq1}),
\begin{eqnarray}\label{con1}
\frac{y_n^{(i)}(I_1,J_1)}{y_n^{(i+1)}(I_1,J_1)}&=&\left(\sum_{k=2}^{\infty}h_k\sum_{j=2}^{k}\big(a_{k,j}-E_{k,j}(\bn{1}-\bm{y}_n(I_1,J_1))_{i\rightarrow i+k-1}\big)\frac{y_n^{(i+j-1)}(I_1,J_1)}{y_n^{(i+1)}(I_1,J_1)}\right)\nonumber\\
& &\boldsymbol{\cdot}\left(1-\sum_{k=1}^{\infty}h_k\big(a_{k,1}-E_{k,1}(\bn{1}-\bm{y}_n(I_1,J_1))_{i\rightarrow i+k-1}\big)\right)^{-1}.
\end{eqnarray}
Since $J_1(n)\to \infty$, there exists $N_1$ such that for $n\ge N_1$,  $J_1(n)\geq N_0$. From Lemma \ref{lem6}, we have $y_n^{(m)}(I_1,J_1)> y_n^{(m+1)}(I_1,J_1)$ for any $m>0$. Observe that \begin{eqnarray}\label{ee} E_{k,i}(\bn{1}-\bm{y}_n(I_1,J_1))_{i\rightarrow i+k-1}\ge 0.\end{eqnarray}
Then
\begin{eqnarray}
\frac{y_n^{(i)}(I_1,J_1)}{y_n^{(i+1)}(I_1,J_1)}\leq\frac{M-M_1}{1-M_1}<\infty,\label{eq3}
\end{eqnarray}
for any $n\geq N_1$ and $1\le i\leq I_1(n)$.

For~$i\geq I_1(n)+1$,
\begin{eqnarray}\label{yy}
\frac{y^{(i)}_n(I_1,J_1)}{y^{(i+1)}_n(I_1,J_1)}=\frac{x^{(J_1(n)+i-I_1(n)-1)}}{x^{(J_1(n)+i-I_1(n))}}\leq\sup_{i}\frac{x^{(i)}}{x^{(i+1)}}<\infty.
\end{eqnarray}
Thus, \begin{eqnarray}\label{yy2}\label{y2} \sup_{n\ge N_1}\sup_{i}\frac{y^{(i)}_n(I_1,J_1)}{y^{(i+1)}_n(I_1,J_1)}<\infty. \end{eqnarray}

On the other hand, if we suppose $y^{(i)}>0$ and $y^{(i+1)}=0$ for some $i\geq1$, since~$\bm{y}_n(I_1,J_1)$ converges to~$\bm{y}$ pointwisely, we then have
\begin{equation}\label{e6.3}
\limsup_{n\to \infty}\frac{y^{(i)}_n(I_1,J_1)}{y^{(i+1)}_n(I_1,J_1)}=\infty,
\end{equation}
which contradicts to (\ref{yy2}). Consequently, either $\bm{y}=\bn{0}$ or $y^{(i)}\neq0$ for any $i$.

From above discussion, to prove $\bm{y}\neq\bn{0}$, we only need to prove $y^{(1)}\neq0$.

$\bullet$ We second prove that if $y^{(1)}=\lim_{n\rightarrow\infty}y_n^{(1)}(I_1, J_1)=0$, then ~$\lim_{n\rightarrow\infty}\parallel\bm{y}_n(I_1,J_1)\parallel_{l_2}=0$. 
 First, from $y_n^{(i+1)}(I_1, J_1)\le y_n^{(i)}(I_1,J_1)$, we have for any $i\ge 1$, 
\begin{equation}\label{e6.300}
y^{(i)}=\lim_{n\rightarrow\infty}y_n^{(i)}(I_1, J_1)=0.
\end{equation}
 Letting $n>N_0$ and
$$\alpha_{n,i}:=\frac{y_n^{(i+1)}(I_1,J_1)}{y_n^{(i)}(I_1,J_1)},$$
then $\alpha_{n,i}\leq1$ follows by Lemma \ref{lem6} and $\bm{x}\in\mathcal{H}(\gamma)$. From (\ref{eq3}) and (\ref{yy}), there exists constant $\alpha_{\text{inf}}<\gamma$ such that $\alpha_{n,i}>\alpha_{\text{inf}}>0$. Let
\[
U_{n,i}:=\sum_{j=1}^{\infty}M_j\prod_{l=1}^{j-1}\alpha_{n,i+l-1}.
\]
For $i\geq I_1(n)+1$, then $y_n^{(i)}(I_1,J_1)=x^{(i-I_1(n)-1+J_1(n))}$. It is easy to prove that
\[
\lim_{n\to\infty}\sup_{i\geq I_1(n)+1}|U_{n,i}-1|=0.
\]

For $1\le i\leq I_1(n)$, by  (\ref{yni}), we have
\begin{eqnarray}\label{con2}
\frac{y_n^{(i)}(I_1,J_1)}{y_n^{(i)}(I_1,J_1)}&=&\sum_{k=1}^{\infty}h_k\sum_{j=1}^{k}\big(a_{k,j}-E_{k,j}(\bn{1}-\bm{y}_n(I_1,J_1))_{i\rightarrow i+k-1}\big)\prod_{l=1}^{j-1}\alpha_{n,i+l-1}=1,
\end{eqnarray}
where $\prod_1^0=1$.  Notice that
\[
U_{n,i}=\sum_{k=1}^{\infty}h_k\sum_{j=1}^{k}a_{k,j}\prod_{l=1}^{j-1}\alpha_{n,i+l-1}.
\]
Combining with (\ref{con2}), we have
\begin{eqnarray}\label{con20}
U_{n,i}-1=\sum_{k=1}^{\infty}h_k\sum_{j=1}^{k} E_{k,j}(\bn{1}-\bm{y}_n(I_1,J_1))_{i\rightarrow i+k-1} \prod_{l=1}^{j-1}\alpha_{n,i+l-1}.
\end{eqnarray}
Meanwhile, $0\le E_{k,j}(\bn{1}-\bm{y}_n(I_1,J_1))_{i\rightarrow i+k-1}\le a_{k,j}$ and $\alpha_{n,i+l-1}\le 1$. Clearly, by Lemma~\ref{lem6} and the definition of $\mathcal{H}(\gamma)$, we have
$$\bn{1}-\bm{y}_n(I_1,J_1))_{i\rightarrow i+k-1}\ge \bn{1}-\bm{y}_n(I_1,J_1))_{1\rightarrow k}.$$
Since $\bm{E}(\bm{s})$ is non-increasing  in $\bm{s}$ (with respect to the partial order induced by ``$\leq$")  and $\bm{E}(\bm{s})\to\bn{0}$ as $\bm{s}\to\bn{1}$. Then by (\ref{e6.300}), for any $i,k,j>0$,
  \begin{eqnarray*}
E_{k,j}(\bn{1}-\bm{y}_n(I_1,J_1))_{i\rightarrow i+k-1}
\le E_{k,j}(\bn{1}-\bm{y}_n(I_1,J_1))_{1\rightarrow k}\to 0,~\mbox{as}~n\longrightarrow\infty.
\end{eqnarray*}
By our assumption {\bf A3},
$$
\sum_{k=1}^{\infty}h_k\sum_{j=1}^{k} a_{k,j}=\sum_{k=1}^{\infty}M_j<\infty.
$$
Applying the dominated convergence theorem in (\ref{con20}), we obtain
$$
\lim_{n\to\infty}\sup_{1\le i\leq I_1(n)}|U_{n,i}-1|=0.
$$
Using Lemma \ref{lem9.5} yields that $\alpha_{n,i}$ converges to $\gamma$ uniformly for $i$.

 Hence, for any $\epsilon>0$ with $\gamma+\epsilon<1$, there exists $N(\epsilon)$ such that for $n>N(\epsilon)$ we have
\[
\sup_{i}\alpha_{n,i}=\sup_{i}\frac{y_n^{(i+1)}(I_1,J_1)}{y_n^{(i)}(I_1,J_1)}<\gamma+\epsilon<1,
\]
which implies
\[
\parallel\bm{y}_n(I_1,J_1)\parallel_{l_2}\leq C_1\cdot y_n^{(1)}(I_1,J_1)
\]
for some constant $C_1$. Therefore $y^{(1)}=\lim_{n\rightarrow\infty}y_n^{(1)}(I_1, J_1)=0$ leads to  $$\lim_{n\rightarrow\infty}\parallel\bm{y}_n(I_1,J_1)\parallel_{l_2}=0.$$

From Lemma~\ref{lem8}, $\lim_{n\rightarrow\infty}\parallel\bm{y}_n(I_1,J_1)\parallel_{l_2}\in(0,\infty)$,  there is a contradiction.  Consequently, we conclude that  $y^{(1)}=\lim_{n\rightarrow\infty}y_n^{(1)}(I_1, J_1)\neq0$ and hence~$y^{(i)}\neq 0$ for all $i\ge 1$. The proof is completed.\qed

\textbf{Proof of Theorem \ref{th1}:}

From Lemma~\ref{lem10}, there exist $\bm{y}\in (0,1]^{\mathbb{N}}$,  $I_1, J_1\in\mathcal{A}$, such that $\bm{y}\neq (1-q)\bn{1}$ and $\bm{y}_n(I_1,J_1)$ converges to~$\bm{y}$ pointwisely. Now, we prove that $\bn{1}-\bm{y}$ is the fixed point of~$\bm{F}(\cdot)$ which is clearly equivalent with $\bm{y}$ is the fixed point of~$\bm{T}(\cdot)$.\par
For any~$i\ge 1$, it holds that
\begin{eqnarray}
& &|T^{(i)}(\bm{y})-y^{(i)}|\nonumber\\
&\leq&|T^{(i)}(\bm{y})-T^{(i)}\big(\bm{y}_n(I_1,J_1)\big)|+|T^{(i)}\big(\bm{y}_n(I_1,J_1)\big)-y^{(i)}_n(I_1,J_1)|+
|y^{(i)}_n(I_1,J_1)-y^{(i)}|\nonumber\\
&=:& K_1+K_2+K_3.\label{e6.4}
\end{eqnarray}
From Lemma~\ref{lem4}, $T(\cdot)$ is pointwisely continuous in $l_2$. Then
\[
\lim_{n\rightarrow\infty}K_1=\lim_{n\rightarrow\infty}|T^{(i)}(\bm{y})-T^{(i)}\big(\bm{y}_n(I_1,J_1)\big)|=0.
\]
From Lemma \ref{lem7},
\[
\lim_{n\rightarrow\infty}K_2=\lim_{n\rightarrow\infty}|T^{(i)}\big(\bm{y}_n(I_1,J_1)\big)-y^{(i)}_n(I_1,J_1)|=0.
\]
It follows from~$\bm{y}_n(I_1,J_1)$ converges to $\bm{y}$ pointwisely that
\[
\lim_{n\rightarrow\infty}K_3=\lim_{n\rightarrow\infty}|y^{(i)}_n(I_1,J_1)-y^{(i)}|=0.
\]
Therefore, letting $n\rightarrow\infty$ in (\ref{e6.4}) yields~$T^{(i)}(\bm{y})=y^{(i)}$ for any $i\ge 1$. Hence $\bm{y}$ is the fixed point of~$\bm{T}(\cdot)$ and  clearly $\bm{y}\in(0,1-q)^{\mathbb{N}}$ follows by $\bm{F}(q\bn{1})=q\bn{1}$, where $q\bn{1}$ is the global extinction probability of $\{\bm{Z}_n;n\geq0\}$.

Next, define $\bm{y}_1=(T_{-1}^{(1)}(\bm{y}),\bm{y})$ and $\bm{y}_i=(T_{-1}^{(1)}(\bm{y}_{i-1}),\bm{y}_{i-1})$ for $i>1$. Clearly from Lemma \ref{lem1} and $\bm{F}(q\bn{1})=q\bn{1}$, $T_{-1}^{(1)}(\bm{y}_{i-1})<1-q$ as long as $y^{(1)}_{i-1}<1-q$. Thus, from the definition of $T_{-1}^{(1)}(\cdot)$ and $\bm{T}(\bm{y})=\bm{y}$,  we obtain that $\bm{y}_i~(i\ge 1)$ are also the fixed points of~$\bm{T}(\cdot)$, $\bm{y}_i\in(0,1-q)^{\mathbb{N}}$ and $\bm{y}_{i+1}>\bm{y}_{i}$ for $i\ge 1$. Then there are at least countably many fixed points of $\bm{T}(\cdot)$ and $\bm{F}(\cdot)$.

 Next, if $\bm{r}$ is a fixed point of $\bm{F}(\cdot)$ and $\bm{r}\notin\Theta$. From \cite[Lemma 3.3]{MJE1}, we know that $\sup_{i}r^{(i)}=1$ and $r^{(i)}\neq1$ for all $i$.  Define $b_i=\frac{1-r^{(i+1)}}{1-r^{(i)}}$ for $i\ge 1$. From (\ref{eq3}), $1>b_i>\frac{1-M_1}{M-M_1}$. Since $\bm{r}$ is the fixed point, by (\ref{eq}) we obtain
\[
1-r^{(i)}=\sum_{k=1}h_k(1-f_k(\bn{1}-\bm{r}_{i\to i+k-1})).
\]
By similar calculation with (\ref{eq2}), dividing $1-r^{(i)}$ in both side yields
\[
1=\sum_{k=1}^{\infty}h_k\sum_{j=1}^{k}\big(a_{k,j}-E_{k,j}(\bn{1}-\bm{r}_{i\rightarrow i+k-1})\big)\prod_{l=1}^{j-1}b_{i+l}.
\]
Hence~$\sum_{k=1}^{\infty}M_k\prod_{l=1}^{k-1}b_{i+l}$ converges to $1$ as $i\to\infty$. Same with the proof of Lemma \ref{lem10}, let
$$
U_{n,i}=\sum_{k=1}^{\infty}M_k\prod_{l=1}^{k-1}b_{n+i+l}.
$$
It satisfies the condition of Lemma \ref{lem9.5} and hence $\lim_{i\to\infty}b_i=\gamma$.
The proof is completed.
\qed

\section*{Acknowledgement}
\par
This work was supported by the National
Natural Science Foundation of China (Grant No. 12271043) and the National Key Research and Development Program of China (No. 2020YFA0712900)

\end{document}